\documentclass[a4paper]{article}
\usepackage{amsmath}
\usepackage{amsbsy}
\usepackage{amsfonts}
\usepackage{amstext}
\usepackage{amssymb}
\usepackage{mathbbol}
\usepackage{wasysym}
\usepackage{eucal}
\usepackage{mathtools}
\usepackage[amsmath,thmmarks]{ntheorem}		
\usepackage{stackrel}
\usepackage{youngtab}

\usepackage{tikz}
   \usetikzlibrary{matrix}
   \usetikzlibrary{calc}
   \usetikzlibrary{decorations,snakes}

\usepackage{float}
\usepackage{subfigure}

\usepackage{enumitem}
\usepackage{mathtools}
\usepackage[a4paper, left=2.4cm, right=2.4cm, top=3.0cm, bottom=2.5cm, bindingoffset=0cm]{geometry}

\usepackage{color}

\definecolor{grau}{rgb}{0.5,0.5,0.5}
\newcommand{\emphdef}[1]{{\textit{\textbf{#1}}}}

\newcommand{\RR}{\ensuremath{{\mathbb R}}}
\newcommand{\QQ}{\ensuremath{{\mathbb Q}}}
\newcommand{\ZZ}{\ensuremath{{\mathbb Z}}}
\newcommand{\CC}{\ensuremath{{\mathbb C}}}

\newcommand{\field}{\ensuremath{{\Bbbk}}}

\newcommand{\Aut}{\ensuremath{\text{Aut}}}
\renewcommand{\ker}{\ensuremath{\text{ker}}}

\newcommand{\liesl}{\ensuremath{{\mathfrak s \mathfrak l}}}

\newcommand{\spann}[1]{\ensuremath{\text{span}_{#1} \:}}
\renewcommand{\dim}{\ensuremath{\text{dim}}}

\newcommand{\aideal}{\ensuremath{{\mathfrak a}}}

\newcommand{\supp}{\ensuremath{\text{Supp}\:}}

\newcommand{\wirk}{\ensuremath{{}_{^{_{\bullet}}}}}

\newcommand{\ua}{\ensuremath{^{\ast}}}

\renewcommand{\i}{\ensuremath{^{-1}}}

\newcommand{\p}{\ensuremath{^{\prime}}}

\renewcommand{\t}[1]{\ensuremath{\tilde{#1}}}

\newcommand{\h}[1]{\ensuremath{\widehat{#1}}}
\newcommand{\ol}[1]{\ensuremath{\overline{#1}}}

\newcommand{\up}{\ensuremath{^{{}^\text{UP}}}}
\newcommand{\low}{\ensuremath{^{{}^\text{LOW}}}}

\newcommand{\lrar}{\leftrightarrow}

\newcommand{\rar}{\rightarrow}

\newcommand{\srar}{\twoheadrightarrow}

\newcommand{\actson}{\ensuremath{\circlearrowright}}

\newcommand{\LieG}{\ensuremath{{\mathfrak g}}}
\newcommand{\Cartan}{\ensuremath{{\mathfrak h}}}

\newcommand{\UE}{\ensuremath{{\cal U}}}

\newcommand{\brac}[1]{\ensuremath{\langle #1 \rangle}}

\newcommand{\mspec}{\ensuremath{\text{mspec}}}			
\newcommand{\maxi}{\ensuremath{{\mathfrak m}}}			
\newcommand{\anni}[1]{\ensuremath{\text{Ann}_{#1}}}	

\theoremstyle{plain}
\theoremseparator  {.} 
\theoremheaderfont {\bfseries}
\theorembodyfont   {\normalfont}
\theoremnumbering  {arabic}
\makeatletter
\newtheoremstyle{normal}%
{\item[\hskip\labelsep \theorem@headerfont ##1\ ##2\theorem@separator]\normalfont}%
{\item[\hskip\labelsep \theorem@headerfont ##1\ ##2]{\theorem@headerfont (##3)}\theorem@separator\ \normalfont}
\newtheoremstyle{nonumber}%
{\item[\theorem@headerfont\hskip\labelsep ##1\theorem@separator]\normalfont}%
{\item[\theorem@headerfont\hskip \labelsep ##1]{\theorem@headerfont (##3)}\theorem@separator\ \normalfont}
\makeatother

\theoremstyle{normal}
\newtheorem{thm}               {Theorem} [subsection]
           
\newtheorem{lemma}      [thm]  {Lemma}

\newtheorem{cor}        [thm]  {Corollary}

\newtheorem{defi}       [thm]  {Definition}
\newtheorem{prop}       [thm]  {Proposition}
\newtheorem{bem}        [thm]  {Remark}
\newtheorem{bsp}        [thm]  {Example}

\theoremstyle{nonumber}
\theoremsymbol{$\square$}
\newtheorem{bew}{Proof}
\theoremsymbol{}
\theoremstyle{nonumber}
\newtheorem{thmnn}{Theorem}

\theoremstyle{nonumber}
\theoremstyle{nonumber}
\theoremstyle{nonumber}
\theoremsymbol{$\smiley$}

\theoremsymbol{}

\begin{document}

\title{Duflo Theorem for a Class of Generalized Weyl Algebras}
\author{Joanna Meinel}
\date{}
\maketitle

\begin{abstract}
For a special class of generalized Weyl algebras, we prove a Duflo theorem stating that the annihilator of any simple module is in fact the annihilator of a simple highest weight module.
\end{abstract}

\section*{Introduction}
Let $\field$ be an algebraically closed field of characteristic $0$. For the universal enveloping algebra of a semisimple Lie algebra over $\field$, Duflo's Theorem \cite{duflo} states that all its primitive ideals (i.e. the annihilators of simple modules) are given by the annihilators of simple highest weight modules. In contrast, the simple modules themselves are far from being classified in general.
Fortunately, for several other classes of algebras the notion of a highest weight module makes sense and the analogue of Duflo's theorem holds:

In \cite{smith}, Smith introduced a family of algebras similar to $\UE(\liesl_2)$. These are $\CC$-algebras generated by three elements $E,F,H$ subject to the relations $[H,E]=E$, $[H,F]=-F$ and $[E,F]=f(H)$ where $f$ can be any polynomial. They share many properties with $\UE(\liesl_2)$ (which is of course included in this family for $f(H)=2H$). In particular it is straightforward to generalize the notion of highest weight modules to these algebras and indeed all primitive ideals are given by annihilators of highest weight modules (see \cite[Theorem 3.3]{smith}).

For classical simple Lie superalgebras, Musson defines $\ZZ_2$-graded highest weight modules depending on a choice of a triangular decomposition. Then all $\ZZ_2$-graded primitive ideals in the universal enveloping algebra of a classical simple Lie superalgebra are given by the annihilators of $\ZZ_2$-graded simple highest weight modules (see \cite[Theorem 2.2]{musson-super}).

In \cite{mvdb}, Musson and Van den Bergh introduce algebras that, roughly speaking, allow a weight space decomposition with weight spaces cyclic over a commutative subalgebra. This class of algebras is closed under taking certain graded subalgebras, tensor products and central quotients. They show that (under some further assumptions, see Theorem \ref{thm:mvdb} for details) all prime, hence all primitive ideals are given by the annihilators of simple weight modules. In particular, this applies to localizations of Weyl algebras and their central subquotients (see \cite[Chapter 6]{mvdb}). Note that for a classical Weyl algebra, given by differential operators on a polynomial ring in $n$ variables, the primitive ideals are not very interesting: These algebras are simple, i.e. the only proper twosided ideal is the zero ideal. Since an annihilator is always twosided, the only primitive ideal of a classical Weyl algebra is the zero ideal.

Now it is natural to ask whether an analogous statement holds for generalized Weyl algebras, a class of algebras that includes many interesting examples, in particular Smith's generalizations of $\UE(\liesl_2)$. These noncommutative algebras are generated by a commutative $\field$-algebra $R$ together with $2n$ elements $X_1,\ldots,X_n$, $Y_1,\ldots,Y_n$. For the relations see Section \ref{sec:prelim}. They are $\ZZ^n$-graded by setting $\deg(X_i)=e_i$, $\deg(Y_i)=-e_i$ where $e_i$ denotes the $i$-th standard basis vector in $\ZZ^n$. Each graded component is a cyclic $R$-module. In this situation, we can define highest weight modules and formulate a Duflo theorem. We prove it for a special class of generalized Weyl algebras using a theorem by \cite{mvdb} that relates the annihilator of a simple weight module to its support and obtain as main result (see Theorem \ref{thm:main}):
\begin{thmnn}
Let $A=R(\sigma,t)$ be a GWA of rank $n$ as defined in Section \ref{sec:prelim} where we assume $R=\field[T_1,\ldots,T_n]$, $\sigma_i(T_j)=T_j-\delta_{ij}b_i$ for  $b_i\in\field\setminus\{0\}$ and   $t_i\in\field[T_i]\subset \field[T_1,\ldots,T_n]$, $t_i\notin\field$.
Then all primitive ideals of $A$, i.e. the annihilator ideals of simple $A$-modules, are given by the annihilators of simple highest weight $A$-modules $L(\maxi)$ of highest weight $\maxi\in\mspec(R)$. 
\end{thmnn}
In Section \ref{sec:prelim} we recall the definition of generalized Weyl algebras, define highest weight modules and discuss graded modules over generalized Weyl algebras. We characterize moreover the highest weight modules as those modules with a locally nilpotent action of the $X_i$. In Section \ref{sec:main} we formulate and prove the main theorem. The principal tool is the Duflo type theorem using \emph{weight} modules from \cite{mvdb}. We show it applies to our situation and improve it by showing that it is enough to consider the much smaller class of \emph{highest weight} modules (as in the classical Duflo theorem). In Section \ref{sec:ex} we finally give some examples to illustrate the relationship between the annihilator and the support of simple highest weight modules.

\subsubsection*{Acknowledgements}
This work is part of my PhD studies at the MPIM Bonn. The subject was suggested to me by Volodymyr Mazorchuk and carried out during my stay at the University of Uppsala which I enjoyed a lot -- I am deeply grateful for numerous interesting discussions and the kind hospitality! I would like to thank my advisor Catharina Stroppel and also Jonas Hartwig for many helpful remarks.
I am supported by a scholarship of the Deutsche Telekom Stiftung.

\section{Generalized Weyl algebras and graded modules}\label{sec:prelim}
\subsection{Generalized Weyl algebras}
Fix a base field $\field=\ol{\field}$ of characteristic $0$. Fix a unital associative commutative $\field$-algebra $R$ that is a noetherian domain. 
Given $n$ nonzero elements $t=(t_1,\ldots,t_n)$ in $R$ and $n$ pairwise commuting algebra automorphisms $\sigma=(\sigma_1,\ldots,\sigma_n)$ in $\Aut(R)$, define the corresponding \emphdef{generalized Weyl algebra (GWA)} $A = R(\sigma, t)$ of rank $n$ as follows: It is the $\field$-algebra generated over $R$ by $2n$ generators $X_i,Y_i$, $1\leq i\leq n$
with relations given by
\begin{align*}
X_i r\ &=\ \sigma_i(r)X_i, \quad &X_i Y_i\ &=\ \sigma_i(t_i), \quad&[X_i,X_j]\ &=\ 0,\\
Y_i r\ &=\ \sigma_i\i(r)Y_i, &Y_i X_i\ &=\ t_i, &[Y_i,Y_j]\ &=\ 0,\\
&&&&[X_i,Y_j]\ &=\ 0
\end{align*}
for all $1\leq i,j\leq n$ with $i\neq j$ and all $r\in R$. It was introduced originally by Bavula in \cite{bav}.
Assume furthermore that $\sigma_i(t_j)=t_j$ for all $i\neq j$, then  
$A=\bigoplus\limits_{\alpha\in\ZZ^n} R\cdot a^\alpha$ is a free left and right $R$-module with generators
$$a^\alpha = a_1^{\alpha_1}\cdot\ldots\cdot a_n^{\alpha_n},\quad a_i^{\alpha_i} = \begin{cases}X_i^{\alpha_i} &\text{ for }\alpha_i\geq 0\\ Y_i^{|\alpha_i|} &\text{ for }\alpha_i< 0,\end{cases}$$
see eg. \cite[Section 1.2]{bavbek} or \cite[Lemma 2.3]{bo}. 
Denote $A_\alpha=R\cdot a^\alpha$. Since $A_\alpha\cdot A_\beta\in A_{\alpha+\beta}$, any GWA $A$ is a $\ZZ^n$-graded algebra with $\deg(X_i)=e_i$ and $\deg(Y_i)=-e_i$ where we denote by $e_i$ the $i$-th standard basis vector of $\ZZ^n$, see eg. \cite[Section 1.1]{bav}. The degree $0$ part of $A$ is given by $A_0=R$.
Notice that the $\sigma_1,\ldots,\sigma_n$ from the defining data of a GWA $A=R(\sigma,t)$ give a $\ZZ^n$-action on $R$ by $e_i\mapsto\sigma_i$.
Write $\sigma^\alpha = \sigma_1^{\alpha_1}\cdot\ldots\cdot\sigma_n^{\alpha_n}$.
\begin{lemma}\label{lem:basics}
Let $A=R(\sigma,t)$ be a GWA of finite rank. Then
$A$ is a left and right noetherian domain, and the tensor product over $R$ or $\field$ of two GWA's is again a GWA.
\end{lemma}
\begin{bew}
$A$ is a domain since $R$ is a domain and $t_i\neq 0$, and noetherianity is a consequence of $R$ being noetherian (see \cite[Proposition 1.3]{bav}, \cite{bo}). The last statement is obvious (see \cite[1.1]{bav}).
\end{bew}

\subsection{A special class of GWA's}\label{sec:setup}
We confine ourselves to the study of GWA's with base ring $R=\field[T_1,\ldots,T_n]$, automorphisms $\sigma_i(T_j)=T_j-\delta_{ij}b_i$ for  $b_i\in\field\setminus\{0\}$ and $t_i\in\field[T_i]\subset \field[T_1,\ldots,T_n]$, $t_i\notin\field$.
This is the tensor product of $n$ GWA's of rank $1$ over the polynomial ring in one variable $\field[T]$, with $\sigma\in\Aut(\field[T])$ of the form $T\mapsto T-b$ for some $b\neq0$ in $\field$ and a nonconstant element $t\in \field[T]$. As $\field$ is algebraically closed, we can factorize $t=(T-z_1)\cdot\ldots\cdot(T-z_s)$ for some $z_1,\ldots,z_s\in \field$ (multiplying $t$ by some nonzero scalar would give an isomorphic GWA, so we can assume this scalar is $1$). 
\begin{bem}
With this choice of $\sigma_1,\ldots,\sigma_n$ the $\ZZ^n$-action on $R$ is free on $R\setminus \field$ (on $\field$ the action is trivial). Additionally, the $\ZZ^n$-action on $\mspec(R)$ given by $\alpha\wirk\maxi:=\sigma^\alpha(\maxi)$ is free. As freeness is defined pointwise, every orbit $\{\sigma^\alpha(\maxi)\ |\ \alpha\in\ZZ\}$ is infinite.
So we only deal with pure translations, i.e. $a=1$ in a general automorphism $\sigma:T\mapsto aT-b$, $a\neq0$, of $\field[T]$. For the application of \cite{mvdb}, we need to work with $\ZZ$-lattices, and we want to keep things easy.
\end{bem}

\subsection{Weight modules}\label{sec:weight}
In this section, $A=R(\sigma,t)$ can be any GWA. By a module, we always mean a left module unless stated otherwise. Denote by $\mspec(R)$ the set of maximal ideals of $R$. For $\maxi\in\mspec(R)$ define the \emphdef{$\maxi$-weight space} of an $A$-module $M$ to be 
$$M_\maxi\ =\ \{v\in M\ |\ \maxi \cdot v =0\}$$
and say that $M$ is a \emphdef{weight module} if $M$ decomposes as vector space into its weight spaces $M=\sum_{\maxi\in\mspec(R)} M_\maxi$. Define the \emphdef{support} of the weight module $M$ to be
$$\supp(M)\ =\ \{\maxi\in\mspec(R)\ |\ M_\maxi\neq0\}.$$
Furthermore, for a weight $A$-module $M$ we have $X_i(M_\maxi)\subset M_{\sigma_i(\maxi)}$ and $Y_i(M_\maxi)\subset M_{\sigma_i\i(\maxi)}$. In other words, $A_\alpha\cdot M_\maxi\subset M_{\sigma^\alpha(\maxi)}$ for $\alpha\in\ZZ^n$.
$M$ is called a \emphdef{highest weight} $A$-module if it is generated as $A$-module by $M_\maxi$ and $X_i\cdot M_\maxi=0$ for all $1\leq i\leq n$. In particular, for the support of a highest weight module $M$ we have $\supp(M)\subset \{\sigma^\alpha(\maxi)\ |\ \alpha\in\ZZ^n_{\leq0}\}$.
\begin{lemma}\label{lem:weight} Let $A$ be a GWA of finite rank.
\begin{enumerate}[label=\roman{*}), ref=(\ref{lem:weight}.\roman{*})]
\item Let $M$ be a weight $A$-module. Then $M=\ \bigoplus\limits_{\mathclap{\maxi\in\mspec(R)}}\ M_\maxi$.
\item Let $M$ be a weight $A$-module, $U\subset M$ some $A$-submodule. Then $U$ and hence $M/U$ inherit the weight decomposition from $M$, i.e. $U$ is a homogeneous submodule.
\item Let $M,N$ be weight $A$-modules and $f: M\rar N$ be a homomorphism of $A$-modules. Then $f(M_\maxi)\subset N_\maxi$.
\end{enumerate}
\end{lemma}
\begin{bew} 
The proof is standard.
\begin{enumerate}[label=\roman{*}), ref=((\roman{*})]
\item\label{lem:weight1} Let $v_1+\ldots+v_n=0$ with $v_i\in M_{\maxi_i}$, i.e. $\maxi_i\cdot v_i=0$, and assume $\maxi_i\neq\maxi_j$ for all $i\neq j$. In particular, $\prod\limits_{i\neq j}\maxi_i\not\subset\maxi_j$ because $\maxi_j$ is maximal and hence prime. Each $v_j$ is zero:
As $\prod\limits_{i\neq j}\maxi_i\ni r$ annihilates all $v_i$, $i\neq j$, we get $0=r\cdot(v_1+\ldots+v_n)=r\cdot v_j$. So $v_j$ is annihilated by $\prod\limits_{i\neq j}\maxi_i$ and $\maxi_j$ which generate the whole $R$. In particular $1\cdot v_j=v_j=0$ and the sum is direct.
\item We have to check that $U=\ \bigoplus\limits_{\mathclap{\maxi\in\mspec(R)}}\ U_\maxi$ with $U_\maxi:=U\cap M_\maxi$. Decompose $v\in U$ as element of $M$ into $v=v_1+\ldots+v_n$ with nonzero $v_j\in M_{\maxi_j}$. We show by a diagonal argument that $v_j\in U$ for all $j$.
Take some element $r:=\prod\limits_{i\neq j}r_i$, where the $r_i$ are some nonzero elements of the maximal ideals $\maxi_i$. Hence $r$ is nonzero, $r\cdot v_j\neq 0$ and $r\notin \maxi_j$. Thus there is some $r\p\in R$ with $r\p r=1\in \field\cong R/\maxi_j$.
We get $r\p r\cdot v=r\p r\cdot v_j =v_j$ because $\maxi_j$ annihilates $v_j$. Therefore $v_j\in U$.
It follows that $M/U$ is isomorphic to $\bigoplus\limits_{\maxi\in\mspec(R)}M_\maxi/U_\maxi$.
\item Since $f$ is an $A$-module homomorphism, $\maxi\cdot f(v)=f(\maxi\cdot v)$ for all $v\in M$. Hence $\maxi\cdot f(M_\maxi)=0$, in other words, $f(M_\maxi)\subset N_\maxi$.
\end{enumerate}
\end{bew}
From the lemma it follows that the weight $A$-modules together with $A$-module homomorphisms that preserve the weight spaces form a full abelian subcategory of the category of left $A$-modules.

\subsection{A characterization of highest weight modules for special GWA's}
Here, $A$ is a special GWA as defined in Section \ref{sec:setup}. The following lemma characterizes highest weight $A$-modules. A similar result for Lie algebras can be found in \cite{mz}.
\begin{prop}
Let $M$ be a simple left $A$-module. The following are equivalent:
\begin{enumerate}[label=\roman{*}), ref=(\roman{*})]
\item\label{prop:hw1} $M$ is a highest weight module.
\item\label{prop:hw2} For all $1\leq i\leq n$, the action of $X_i$ on $M$ is locally nilpotent, i.e. for every $v\in M$ there exists a natural number $k_i$ such that $X_i^{k_i}\cdot v=0$.
\item\label{prop:hw3} There exists $v\in M$ such that $X_i$ acts nilpotently on $v$ for all $1\leq i\leq n$.
\end{enumerate}
\end{prop}
\begin{bew}
\ref{prop:hw1}$\Rightarrow$\ref{prop:hw2}: Let $M$ be a highest weight module with highest weight $\maxi$ and weight space decomposition $M=\bigoplus\limits_{\alpha\in\ZZ^n_{\leq0}}M_{\sigma^\alpha(\maxi)}$. So any $v\in M$ decomposes as $v=v_{\alpha(1)}+\ldots+v_{\alpha(r)}$ for weight vectors $v_{\alpha(j)}\in M_{\sigma^{\alpha(j)}(\maxi)}$.
In particular, $X^{-\alpha(j)}\cdot v_{\alpha(j)}=a^{-\alpha(j)}\cdot v_{\alpha(j)}\in M_\maxi$, which is the highest weight space, hence
$X_i^{-\alpha(j)_i+1}\cdot v_{\alpha(j)}=0$. Now choose $k_i\in\ZZ$ such that $k_i \geq -\alpha(j)_i+1$ for all $j$. Then
$X_i^{k_i}\cdot v_{\alpha(j)} =0$ for all $j$ and therefore $X_i^{k_i}\cdot v=0$.

\ref{prop:hw2}$\Rightarrow$\ref{prop:hw3}: Clear.

\ref{prop:hw3}$\Rightarrow$\ref{prop:hw1}: Assume we have an element $v\in M$ such that $X_i^{k_i}\cdot v=0$ for some natural numbers $k_i$. We construct a nonzero element $v\p$ in $M$ that is annihilated by all $X_i$ and a maximal ideal $\maxi$.
Since $M$ is simple, this suffices to prove that 
$$M\ =\ A\cdot v\p\ =\ \bigoplus\limits_{\alpha\in\ZZ_{\leq0}^n}A_\alpha v\p\ =\ \bigoplus\limits_{\alpha\in\ZZ_{\leq0}^n} M_{\sigma^\alpha(\maxi)}$$ 
is a highest weight module of highest weight $\maxi$.
Since the $X_i$ commute, we can find for all $i$ a natural number $\beta_i$ such that $0\leq\beta_i<k_i$ and $\t{v}:= X_1^{\beta_1}\ldots X_n^{\beta_n}\cdot v\neq 0$ but $X_i\cdot \t{v}=0$ for all $i$. Hence
$$t_i\cdot \t{v}\ =\ Y_iX_i\cdot \t{v}\ =\ 0.$$
Now according to our assumption $t_i\in \field[T_i]$ is a polynomial, say $t_i=(T_i-a(i)_1)\ldots(T_i-a(i)_{s(i)})$ for some $s(i)\in\ZZ_{>0}$ and $a(i)_r\in\CC$.
So there is a linear factor $(T_i-a(i)_{r(i)})$ such that 
\begin{align*}
v(i)\ :=\ (T_i-a(i)_{r(i)+1})\ldots (T_i-a(i)_{s(i)}) \t{v}\ &\neq\ 0, \quad\text{and}\\ 
(T_i-a(i)_{r(i)})(T_i-a(i)_{r(i)+1})\ldots (T_i-a(i)_{s(i)}) \t{v}\ &=\  0.
\end{align*}
In this way we construct successively nonzero elements $v(1),\ldots,v(n)$ in $M$ that are annihilated by all $X_i$ since they differ from $\t{v}$ only by multiplication with elements in the base ring $R$. Furthermore, $v\p:=v(n)$ is annihilated by the maximal ideal $\maxi:=\left(T_1-a(1)_{r(1)},\ldots,T_n-a(n)_{r(n)}\right)$.
\end{bew}

\subsection{Side remark on generalized gradings}\label{sec:grad}

Theorem \ref{thm:mvdb} describes primitive ideals of graded algebras in terms of annihilators of \emph{graded} simple modules. 
Although GWA's are $\ZZ^n$-graded, their weight modules are not $\ZZ^n$-graded in general. Instead, weight modules $M$ decompose into weight spaces $M_\maxi$ indexed by $\mspec(R)$. It makes sense to think of a weight module as a graded module, but instead of the usual notion of graded modules over a graded algebra, where both objects are graded over the same additive group, one needs to generalize it as follows:
\begin{defi}
Let $G$ be an abelian group and $X$ be a set with $G$-action.
Let $A=\bigoplus\limits_{g\in G}A_g$ be a $G$-graded algebra. Then a $(G\actson X)$-graded module (or a module with $X$-grading respecting the $G$-action) is an $A$-module $M$ with a decomposition $M=\bigoplus\limits_{x\in X}M_x$ such that $A_g\cdot M_x\subset M_{gx}$.
\end{defi}
This kind of graded modules was studied in \cite{NRvanOy}, motivated by $G$-graded modules over the group algebra $\field[G]$ of a group $G$: Take a $\field[G]$-module graded by the group $G$ itself, but consider it now as $\field[H]$-modules for a subgroup $H\subset G$. As a $\field[H]$-module, it is then naturally $(H\actson G)$-graded. In \cite{BeaDas} an equivalence of the category of $(G\actson X)$-graded modules with the module category over a smash product ring is given.

Weight modules over a GWA $A$ are naturally $(\ZZ^n\actson\mspec(R))$-graded because $A_\alpha\cdot M_\maxi\subset M_{\sigma^\alpha(\maxi)}$.
Nevertheless, for our special GWA's it is enough to change the indexing set of both the GWA $A$ and the module $M$ to find a common index set with group structure, with respect to which $M$ is a classically graded $A$-module, see Section \ref{sec:gradinggwa}. So we will work with the classical grading.

\section{Description of weight modules in terms of breaks}

\subsection{Grading of weight modules}\label{sec:gradinggwa}

Let $A$ be again a special GWA as introduced in Section \ref{sec:setup}. Consider the left $A$-module $M(\maxi)=A/A\maxi$. 
As $R$-module it decomposes into 
$$M(\maxi)=\bigoplus\limits_{\alpha\in\ZZ^n} A_\alpha/A_\alpha\maxi,$$
and this decomposition is already a weight space decomposition because
\begin{align*}
A_\alpha/A_\alpha\maxi\ &\cong\ \{m\in A\ |\ \sigma^\alpha(\maxi)\cdot m\in A\maxi\}\ 
\cong\ \{m\in M(\maxi)\ |\ \sigma^\alpha(\maxi)\cdot m =0\}\ 
=\ M(\maxi)_{\sigma^\alpha(\maxi)}.
\end{align*}
Here we use $\sigma^\alpha(\maxi)\cdot A_\alpha = A_\alpha\cdot\maxi$ and that $A$ is graded, so that one can study whether $\sigma^\alpha(\maxi)\cdot m$ is an element of $A\maxi$ for homogeneous $m$. 
For $\maxi=\maxi_a$ and the shorthand notation $M(\maxi_a)_{a\p}=M(\maxi_a)_{\maxi_{a\p}}$ and $\alpha\cdot\beta$ defined componentwise by $(\alpha\cdot b)_i=\alpha_i\cdot b_i$, we obtain
$$M(\maxi_a)\ =\ \bigoplus\limits_{\alpha\in\ZZ^n}M(\maxi_a)_{\sigma^\alpha(\maxi_a)}\ =\ \bigoplus\limits_{\alpha\in\ZZ^n}M(\maxi_a)_{a+\alpha\cdot b}$$
(notice that indeed $\sigma_i(\maxi_a)=\maxi_{a+b_i}$).
This weight space decomposition turns $M(\maxi_a)$ into a graded $A$-module, but only after reindexing the grading of $A$:
The decomposition of $M(\maxi_a)$ does not respect the usual $\ZZ^n$-grading of $A=\bigoplus\limits_{\alpha\in\ZZ^n} A_\alpha$
because $A_\alpha\cdot M(\maxi_a)_{a\p}$ is a subset of $M(\maxi_a)_{a\p+\alpha\cdot b}$ instead of $M(\maxi_a)_{a\p+\alpha}$.
We have to interpret the abstract $\ZZ^n$-grading of the GWA $A$ as a $\ZZ^n\cdot b$-grading coming from the adjoint $R$-action as in \ref{mvdb1}, where we write $\ZZ^n\cdot b=\{\alpha\cdot b\ |\ \alpha\in\ZZ^n\}$. Observe that
\begin{align*}
A_\alpha\ &=\ R\cdot a^\alpha\\ 
&=\ \{a\in A\ |\ r\cdot a = a\cdot\sigma^{-\alpha}(r)\quad \forall\ r\in R\}\\
&=\ \{a\in A\ |\ T_i\cdot a=a\cdot (T_i+\alpha_i b_i)\quad \forall\ 1\leq i\leq n\}\\
&=\ \{a\in A\ |\ [T_i,a]=\alpha_i b_i \cdot a\quad \forall\ 1\leq i\leq n\}\\
&=\ A_{\alpha\cdot b}\qquad\text{in the sense of \ref{mvdb1}}.
\end{align*}
Then thanks to $A_{\alpha\cdot b}\cdot A_{\alpha\p\cdot b}\subset A_{(\alpha+\alpha\p)\cdot b}$,
$$A\ =\ \bigoplus\limits_{\tau\in\field^n} A_\tau\quad\text{with}\quad A_\tau=0 \text{ for }\tau\neq \alpha\cdot b$$
is a $\field^n$-grading of $A$.
Of course we do not change the decomposition of $A$, we only choose a concrete realization for the abstract $\ZZ^n$-grading and added some $0$-summands to $A$. 
With respect to this new grading $M(\maxi)$ is a $\field^n$-graded $A$-module.

Let us recall some further properties of $M(\maxi)$:
\begin{itemize}
\item 
Since $\ol{1}\in A_0/A_0\maxi = M(\maxi)_\maxi$, it follows that $A_\alpha\cdot M(\maxi)_\maxi\ =\ M(\maxi)_{\sigma^\alpha(\maxi)}$ for all $\alpha\in\supp(A)$, therefore the support of $M(\maxi)$ is given by 
$$\supp(M(\maxi_a))\ = \ a+\supp(A)\ =\ a+\ZZ^n\cdot b\ =\ \{a+\alpha\cdot b\ |\ \alpha\in\ZZ^n\},$$
i.e. as subset of $\mspec(R)$, the support equals the whole orbit $\{\sigma^\alpha(\maxi)\ |\ \alpha\in\ZZ\}$.
\item 
Every weight space of $M(\maxi)$ is one-dimensional since
$$M(\maxi)_{\sigma^\alpha(\maxi)}= A_\alpha/A_\alpha\maxi\cong R/\sigma^{-\alpha}(\maxi)\cong\field$$
with $A_\alpha\maxi=R\cdot a_\alpha\cdot\maxi=\sigma^{-\alpha}(\maxi)A_\alpha$.
\item 
Every submodule of $M(\maxi)$ is homogeneous (see Lemma \ref{lem:weight}).
\item 
$M(\maxi)$ has a unique simple top, denoted by $L(\maxi)$. It inherits the grading of $M(\maxi)$. 
Its support is denoted by $\brac{\maxi}:=\supp(L(\maxi))$. We usually consider $\brac{\maxi}$ as subset of $\field^n$.
\item
Notice that although the modules $M(\maxi)$ seem to be very similar (as $\field$-vector spaces they are all isomorphic to $\bigoplus\limits_{\alpha\in\ZZ^n} \field$), two modules $M(\maxi)$, $M(\maxi\p)$ are only isomorphic iff their simple tops $L(\maxi)$ and $L(\maxi\p)$ are isomorphic, too, and the latter are isomorphic iff they have the same support. 
\end{itemize}
The weight space structure of the module $M(\maxi)=A/A\maxi$ and the existence of its simple top were discussed in \cite{bav}.

\subsection{Breaks and the submodule lemma}\label{sec:breaks}
Now recall how the submodules of $M(\maxi_a)$ can be described in terms of its support and the breaks therein, see \cite{bav} and \cite{dgo}. Later on we will see how this carries over to the primitive ideals.
\begin{defi}
A maximal ideal $\maxi\in\mspec(R)$ is called a break ideal in direction $i$ if $t_i\in\maxi$. 
\end{defi}
It deserves this name since the module $M(\maxi_a)$ `breaks' into submodules precisely between its weight spaces $M(\maxi_a)_\maxi$ and $M(\maxi_a)_{\sigma_i(\maxi)}$ for break ideals $\maxi$:
\begin{lemma} \label{lem:break}
Let $M=M(\maxi_a)$ for some $\maxi_a\in\mspec(R)$. Let $\maxi$ be in the support of $M$. 
If $t_i\in\maxi$ then $X_i$ or $Y_i$ act by $0$ between $M_\maxi$ and $M_{\sigma_i(\maxi)}$. Otherwise, $X_i$ and $Y_i$ act up to nonzero scalars as mutually inverse bijections between the weight spaces.
\end{lemma}
\begin{bew} 
Every weight space of $M(\maxi_a)$ is of the form $M(\maxi_a)_{a+\alpha\cdot b}$ of $M(\maxi_a)$ and in particular one-dimensional. Therefore, 
\begin{align*}
X_i\cdot M(\maxi_a)_{a+\alpha\cdot b}\ =\ \begin{cases}0,\quad &\text{iff }\quad X_i\cdot a^\alpha \in A_{\alpha+e_i}\maxi_a;\\
M(\maxi_a)_{a+(\alpha+e_i)\cdot b},&\text{else.}
\end{cases}
\end{align*}
For $\alpha_i\geq0$, we have $X_i a^\alpha =a^{\alpha+e_i} \notin A_{\alpha+e_i}\maxi_a$.
For $\alpha_i<0$, the defining relations of a GWA give $X_i a^\alpha =\sigma_i(t_i) a^{\alpha+e_i}$. So $X_i a^\alpha \in A_{\alpha+e_i}\maxi_a = \sigma^{\alpha+e_i}(\maxi_a)A_{\alpha+e_i}$ iff $\sigma_i(t_i)\in \sigma^{\alpha+e_i}(\maxi_a)$ (use that $A_\alpha$ is a free $R$-module), iff $t_i\in \sigma^{\alpha}(\maxi_a)$. Similarly, we obtain for $Y_i$ that $Y_i\cdot M(\maxi_a)_{a+\alpha\cdot b}=0$ iff $\alpha_i> 0$ and $t_i\in \sigma^{\alpha-e_i}(\maxi_a)$.
In other words: 
\begin{align*}
X_i\text{ acts by zero on }&M(\maxi)_{\sigma^{\alpha}(\maxi)}\text{ iff }\alpha_i <0 \text{ and } t_i\in \sigma^{\alpha}(\maxi),\\
Y_i\text{ acts by zero on }&M(\maxi)_{\sigma^{\alpha}(\maxi)}\text{ iff }\alpha_i>0  \text{ and } t_i\in \sigma^{\alpha-e_i}(\maxi).
\end{align*}
Together this proves the claim.
\end{bew}
Since $\sigma_j(t_i)=t_i$ for $i\neq j$, a maximal ideal $\maxi$ is a break ideal in direction $i$ iff so is $\sigma_j(\maxi)$. The break ideals that are in the same $\sigma_j$-orbits for $j\neq i$ lie on a common hyperplane.
\begin{defi}
We call a hyperplane in $\field^n$ containing all $\sigma_j$-orbits of $\maxi$ for $j\neq i$ a break in direction $i$.
\end{defi}
Notice that every point inside a break is indeed a break ideal. If we identify once more $\field^n$ with $\mspec(R)$, the breaks correspond to hyperplanes parallel to the coordinate hyperplanes. From Lemma \ref{lem:break} we know that breaks should be interpreted as `forward breaks'. Examples will be given in Section \ref{sec:ex}. 
\begin{lemma}\label{lem:length}
The module $M(\maxi)$ has at most
$2^{\prod\limits_{i=1}^n (1+\text{number of zeroes of }t_i)}$ submodules. The subquotients occur with multiplicity $1$. In particular, $M(\maxi)$ has finite length bounded by $\prod\limits_{i=1}^n (1+\text{number of zeroes of }t_i)$, independent of $\maxi$.
\end{lemma}
\begin{bew}
Every submodule $N$ inherits the weight space decomposition from $M(\maxi)$, and because every weight space of $M(\maxi)$ is at most one-dimensional, 
we have
$$N\ =\ \bigoplus\limits_{\maxi\p\in\supp(N)}M(\maxi)_{\maxi\p}.$$
The submodules are therefore completely determined by their supports, in the sense that $N=N\p$ iff $\supp(N)=\supp(N\p)$.
From the discussion of the breaks we know that $X_i$ and $Y_i$ act as mutually inverse (up to multiplication by elements in $R$) bijections between the weight spaces, unless we encounter a weight space that belongs to a break. If a weight between two successive breaks belongs to the support of $N$, all the other weights between these two breaks do as well. The choice of a submodule is thus equivalent to the choice of the breaks (or no breaks at all) for each coordinate direction $i$. The polynomial $t_i$ is contained in the maximal ideal $\maxi_a=(T_1-a_1,\ldots,T_n-a_n)$ iff $a_i$ is a zero of $t_i$.
In particular, $t_i$ can only be contained in finitely many maximal ideals in the orbit $\{\sigma_i^{\alpha_i}(\maxi)\ |\ i\in\ZZ\}$. 
So there are only finitely many breaks in each direction $i$, and they all occur at zeros of $t_i$.
Since there are $\#(\text{zeroes of }t_i)$ breaks in the $i$-th coordinate direction, the statement of the lemma follows.
\end{bew}
The breaks provide in particular a description of the support of the simple modules $L(\maxi_a)$: We have $\supp(L(\maxi_a))\subset\supp(M(\maxi_a))$. In other words, $\brac{\maxi_a}\subset a+\ZZ^n\cdot b$, i.e. the support consists of lattice points. Since $\ol{1}\in L(\maxi_a)_a$, we know that $a\in\brac{\maxi_a}$. 
Again, $X_i$ and $Y_i$ act (up to multiplication by elements in $R$) as mutually inverse bijections between the weight spaces, unless we encounter a weight space that belongs to a break ideal. So informally speaking $\brac{\maxi_a}$ is given by those weights that can be reached from $a$ without crossing a break. 

More precisely: For every $i$, pick the largest $\gamma_i\low<0$ with $t_i\in \sigma_i^{\gamma_i\low}(\maxi_a)$ and the smallest $\gamma_i\up>0$ with $t_i\in \sigma_i^{\gamma_i\up-e_i}(\maxi_a)$ (if they exist).
Under the isomorphism $\mspec(R)\cong \field^n$, denote the $i$-th coordinate of the image of $\sigma_i^{\gamma_i\low}(\maxi_a)$ by $g_i\low$ and the image of $\sigma_i^{\gamma_i\up-e_i}(\maxi_a)$ by $g_i\up$ (and set $g_i\low=-\infty$ resp. $g_i\up=\infty$ in case this does not exist). Then as a subset of $\field^n$,
\begin{align*}
\brac{\maxi_a}\ =\ \supp(L(\maxi_a))\ &=\ (a+\ZZ^n\cdot b)\ \cap\ \{x\in \field^n\ |\ g_i\low < x_i \leq g_i\up\ \text{for all }i\}\\
&=\ (a+\ZZ^n\cdot b)\ \cap\ \{x\in \field^n\ |\ g_i\low+b_i \leq x_i \leq g_i\up\ \text{for all }i\}.
\end{align*}
As these inequalities involve only one coordinate each, the support has the shape of a rectangle with sides consisting of hyperplanes parallel to the coordinate hyperplanes, in case there exist $g_i\up$ and $g_i\low$ (otherwise drop the corresponding hyperplane from the picture).
Of course $g_i\low$, $g_i\up$ are just two zeroes of $t_i$ chosen such that $g_i\low < a_i \leq g_i\up$ and there is no other zero of the polynomial $t_i$ between them in the lattice $a_i+\ZZ\cdot b_i$. The choice of these zeroes depends on $a$ (so we should really write ${}^a g_i\up$ if it wasn't too much index notation).

\section{Primitive ideals of generalized Weyl algebras}\label{sec:main}

\subsection{The main result}

Let $A$ be the special GWA described in Section \ref{sec:setup}.
Denote by $\anni{A}(M):=\{a\in A\ |\ a\cdot M=0\}$ the annihilator of $M$. It is a twosided ideal of $A$. For a simple $A$-module $L$, the annihilator $\anni{A}(L)$ is called a primitive ideal. Then our main result reads as follows:
\begin{thm}\label{thm:main}
Let $A$ be the GWA of rank $n$ given by $R=\field[T_1,\ldots,T_n]$, $\sigma_i(T_j)=T_j-\delta_{ij}b_i$ for some $b_i\in\field\setminus\{0\}$ and some  $t_i\in\field[T_i]\subset \field[T_1,\ldots,T_n]$, $t_i\notin\field$.
Then all primitive ideals of $A$ are of the form $\anni{A}(L(\maxi))$ for some simple highest weight $A$-module $L(\maxi)$ of highest weight $\maxi\in\mspec(R)$. In other words, there is a bijection
$$\{\anni{A}(L(\maxi))\ |\ \maxi\in\mspec(R)\text{ st. }L(\maxi)\text{ is a highest weight module}\} \ \lrar \ 
\{\text{ primitive ideals of }A\}.$$
\end{thm}
This theorem is analogous to the classical Duflo theorem from \cite{duflo} for the universal enveloping algebra $\UE(\LieG)$ of a semisimple Lie algebra $\LieG$, stating that its primitive ideals are given by the annihilators of highest weight modules $L(\lambda)$ where $\lambda\in\Cartan\ua$ for a Cartan subalgebra $\Cartan\subset\LieG$.
The proof is an application of Theorem \ref{thm:mvdb} from \cite{mvdb}, which we recall in Section \ref{sec:mvdb}. In Section~\ref{sec:gradinggwa} we will give more details about the simple highest weight module $L(\maxi)$. The proof itself follows in Sections \ref{sec:proof} and \ref{sec:refine}.
From the proof it follows that
\begin{cor}
$A$ as above has only finitely many different primitive ideals.
\end{cor}
We give some important examples of the class of GWA's to which Theorem \ref{thm:main} applies:
\begin{bsp}
\begin{enumerate}[label=\roman{*})]
\item Classical Weyl algebras $A_n=\field[x_1,\ldots,x_n,\partial_1,\ldots,\partial_n]$ (see \cite[Example 1.2.(1)]{bav}). Since these algebras are simple, every primitive ideal is zero.
\item The universal enveloping algebra $\UE(\liesl_2) =\CC\langle e,f,h\rangle/\ ([h,e]=2e,\ [h,f]=-2f,\ [e,f]=h)$ is not included in this class of algebras: It is isomorphic to the GWA $\CC[C,H](\sigma,t)$ with $\sigma(H)=H-2$, $\sigma(C)=C$ and $t=\frac14(C-H(H+2))$. The isomorphism is given by $X\mapsto e$, $Y\mapsto f$, $H\mapsto h$ and $C\mapsto c$ where $c=h(h+2)+4fe$ denotes the Casimir element in the universal enveloping algebra. Hence $t$ is mapped to $fe$. 
However, every simple $\liesl_2$-module $L$ has central character, so for every simple module $L$ there is some $\chi\in \CC$ such that $c\cdot v=\chi\cdot v$ for all $v\in L$. Hence we have
\begin{align*}
\{\text{primitive ideals of }\UE(\liesl_2)\}\ &=\ \bigcup\limits_{\chi\in\CC}\{\text{primitive ideals of }\UE(\liesl_2)\text{ that contain }(c-\chi)\}\\
&\lrar\ \bigcup\limits_{\chi\in\CC}\{\text{primitive ideals of }\UE(\liesl_2)/(c-\chi)\}.
\end{align*}
But the central quotient $\UE(\liesl_2)/(c-\chi)$ is isomorphic to the GWA 
$\CC[H](\sigma:H\mapsto H+2,\ t=\frac14(\chi-H(H+2))$ (see \cite[Example 1.2.(3)]{bav}), to which our theorem applies. Hence we recover the Duflo theorem in this case.
\item More generally, for all $k$-algebras $A$ with $\dim(A)<|k|$ it is true that every simple module has central character, see the argument in \cite[Theorem 4.7]{maz-sl2} (it is shown there that the Casimir element $C$ of $\UE(\liesl_2)$ acts by a scalar on any simple $\UE(\liesl_2)$-module, but one can apply exactly the same argument for a central element $C\in A$ of any algebra with $\dim(A)<|\field|$, eg. $A$ with countable dimension and $\field$ uncountable). 
To obtain a Duflo statement for $A$, it is enough to establish a Duflo theorem for all central quotients $A/(Z-\chi(Z))$, where $Z$ denotes the center of $A$ and $\chi\in Z\ua$ is a central character - similarly to the $\UE(\liesl_2)$-example. The primitive ideals in $A/(Z-\chi(Z))$ can then be lifted to ideals in $A$, which are indeed primitive and exactly those primitive ideals of $A$ that contain $(Z-\chi(Z))$ (all simple $A$-modules with central character $\chi$ are lifts of the simple $A/(Z-\chi(Z))$-modules). But notice that if some $X_i$ is central, a simple highest weight $A/(Z-\chi(Z))$-module need not be \emph{highest} weight as $A$-module in the  sense of the defintion given in Section \ref{sec:weight}. However, it seems to be adequate to adapt the notion of a highest weight module so that a central $X_i$ is not supposed to act by $0$ on the `highest weight space'.
\item More generally, Smith's generalizations of $\UE(\liesl_2)$, defined in \cite{smith}, have central quotients that are GWA's in the special class we consider here. The realization as GWA is given in \cite[Example 1.2.(4)]{bav}. The primitive ideals were already described in \cite[Section 3]{smith}.
\item The class of GWA's and all examples discussed in \cite[Section 1.2]{bav}: They agree with our special GWA's, except that the automorphism $\sigma$ is given by translation by $1$ instead of any nonzero $b$. In \cite[Theorem 3.2, 3.8]{bav}, a classification of simple modules for these algebras is given.
\end{enumerate}
\end{bsp}
We confine ourselves to the special class of GWA's because we want the following properties to hold, mainly for the application of Theorem \ref{thm:mvdb}. Some of them could be weakened slightly, but without greater insight and to the cost of additional technical considerations (as illustrated in the enveloping algebra example).
\begin{itemize}
\item The base ring $R$ is in particular noetherian, hence by Lemma \ref{lem:basics} the GWA $A$ is noetherian, too. This is a requirement of Theorem \ref{thm:mvdb}. 
\item The base ring is the polynomial ring and not just a quotient of such since  otherwise we cannot guarantee that there are only finitely many `breaks', see Section \ref{sec:breaks}. But such a finiteness condition is needed in Theorem \ref{thm:mvdb}.
\item To satisfy $\sigma_i(t_j)=t_j$ for $i\neq j$, it is convenient to consider only tensor products of rank $1$ GWA's. 
\item The application of Theorem \ref{thm:mvdb} is only possible for a GWA where $\ZZ^n$ acts freely on $R$, i.e. $\sigma^\alpha=\sigma^\beta$ iff $\alpha=\beta$: This ensures that the graded components $A_\alpha$ are cyclic over $R$, see \ref{mvdb2} below.
\item The grading should come from a weight space decomposition with respect to the adjoint action of $R$ on $A$. In this case, any twosided ideal inherits the grading of $A$, and this is fundamental for Theorem \ref{thm:mvdb}. Therefore in the rank $1$ case, some automorphism of the polynomial ring $\sigma:T\mapsto aT-b$ must be of the form $\sigma:T\mapsto T-b$.
\item Furthermore, $b_i\neq 0$ because otherwise $\sigma_i$ would be trivial. This contradicts the free $\ZZ^n$-action on $R$.
\end{itemize}

\subsection{The result of \cite{mvdb}}\label{sec:mvdb}

We would like to apply the following result of \cite[Theorem 3.2.4]{mvdb}, slightly reformulated:
\begin{thm}\label{thm:mvdb}
Let $\field$ be an algebraically closed field of characteristic $0$. Let $A$ be any unital associative $\field$-algebra satisfying the following assumptions:
\begin{enumerate}[label=\bf{(A\arabic{*})}, ref=(A\arabic{*})]\setcounter{enumi}{0} 
\item\label{mvdb1} $A$ carries a grading $\bigoplus\limits_{\tau\in\field^n} A_\tau$ with $A_{0}= R:= \field[T_1,\ldots,T_n]$ commutative, where the grading comes from the weight space decomposition of $A$ with respect to the adjoint action of $\spann{\field}\{T_1,\ldots,T_n\}$,
$$A_\tau = \{a\in A\ |\ [T_i,a]=\tau_i a\}.$$
\item\label{mvdb2} $R\srar A_\tau = R\cdot a_\tau$ for all $\tau$, i.e. each nonzero $A_\tau$ is generated by one element over $R$.
\item\label{mvdb3} $A$ is graded left noetherian.
\item\label{mvdb4} For a maximal ideal $\maxi\subset R$, the $A$-module $M(\maxi):=A/A\maxi$ has uniformly bounded length, independent of $\maxi$.
\item\label{mvdb5} The number of different Zariski closed sets $\ol{\brac{\maxi}}\subset\field^n$ is finite.

Here, the set $\brac{\maxi}$ is defined as follows: For an algebra $A$ satisfying \ref{mvdb1} and \ref{mvdb2}, the $A$-module $M(\maxi)$ has a weight space decomposition which turns it into a $\field^n$-graded module with $M(\maxi)_{a}:=M(\maxi)_{\maxi_a}$ and $\maxi_a=(T_1-a_1,\ldots,T_n-a_n)$ is the maximal ideal corresponding to $a=(a_1\ldots,a_n)\in\field^n$: Indeed $A_\tau\cdot M(\maxi)_{\alpha}\subset M(\maxi)_{\alpha+\tau}$.
It is easy to see that $M(\maxi)$ has a unique maximal submodule, because a submodule is proper iff it does not contain $\ol{1}\in A/A\maxi$. Hence $M(\maxi)$ has simple top, denoted $L(\maxi)$. It inherits the grading of $M(\maxi)$. Its support is denoted by $\brac{\maxi}:=\supp(L(\maxi))$. We usually consider $\brac{\maxi}$ as subset of $\field^n$.
\item\label{mvdb6} For all $\maxi_\alpha\in\mspec(R)$ and all $\tau\in\supp(A)$ we have $\ol{(\tau+\brac{\maxi})}\cap\ol{\brac{\maxi}} = \ol{(\tau+\brac{\maxi})\cap\brac{\maxi}}$.
\end{enumerate}
Then all prime ideals, hence all primitive ideals of $A$ are of the form $\anni{A}(L(\maxi))=: J(\maxi)$ for some $\maxi\in\mspec(R)$, and
$$\{\ol{\brac{\maxi}}\ |\ \maxi\in\mspec(R)\} \ \lrar \ 
\{J(\maxi)\ |\ \maxi\in\mspec(R)\} \ \lrar \ 
\{\text{primitive ideals of }A\}.$$
The first bijection is given by $J(\maxi)=A\cdot I(\ol{\brac{\maxi}})\cdot A$ where $I(\ol{\brac{\maxi}})=\bigcap\limits_{\maxi\p\in\ol{\brac{\maxi}}}\maxi\p$.
\end{thm}

The formulation of the theorem is slightly modified: In \cite{mvdb} the subalgebra $R$ can be any finitely generated commutative subalgebra.
We will obtain a slight refinement, by finding the above correspondence for \emph{highest} weight modules $L(\maxi)$.

As mentioned in Section \ref{sec:gradinggwa}, the weight space structure of the module $M(\maxi)=A/A\maxi$ and the existence of its simple top were treated for GWA's already in \cite{bav}. But in fact they are a general consequence of conditions \ref{mvdb1} and \ref{mvdb2} (see \cite[Proposition 3.1.7]{mvdb}).

\subsection{The proof of Theorem \ref{thm:main}: Reduction to weight modules}\label{sec:proof}
We now check the conditions of Theorem \ref{thm:mvdb}.

Condition \ref{mvdb1} is valid for any GWA (here we have to use the unusual grading as described in Section \ref{sec:gradinggwa}).

Condition \ref{mvdb2} holds for any GWA with free $\ZZ^n$-action on $\Aut(R)$. For $\sigma_i$ given by translations in coordinate direction $i$, it follows from $\sigma^\alpha=\sigma^\beta$ that $\alpha=\beta$, so the $\ZZ^n$-action on $\Aut(R)$ is indeed free.

Condition \ref{mvdb3} holds for any GWA whose ground ring $R$ is noetherian (Lemma \ref{lem:basics}), in particular in our case where $R=\field[T_1,\ldots,T_n]$ is the polynomial ring.
 
Condition \ref{mvdb4} is satisfied according to Lemma \ref{lem:length}, and the length is uniformly bounded by $\prod\limits_{i=1}^n (1+\text{number of zeroes of }t_i)$.

For the verification of \ref{mvdb5} and \ref{mvdb6}, we first notice that there are only finitely many breaks (i.e. hyperplanes consisting of those points in $\field^n$ that correspond to maximal ideals $\maxi\subset \field[T_1,\ldots,T_n]$ containing one of the $t_i$).
\begin{bem}
In case $\maxi$ is contained in an orbit without breaks, the support of $L(\maxi)$ is the whole orbit $\brac{\maxi}=\supp(A)\cdot\maxi$. For our special choice of GWA's $A$ we have $\supp(A)\cdot\maxi=\ZZ^n\cdot\maxi$ which is dense in $\mspec(R)$, and therefore $\ol{\brac{\maxi}}=\mspec(R)$. So these closures give all the same contribution when we count the different closures to verify \ref{mvdb5}. Also, $\sigma^\alpha(\brac{\maxi})=\brac{\maxi}$ for any $\sigma^\alpha\in\supp(A)$ and so \ref{mvdb6} is satisfied for those $\maxi$.
\end{bem}
For $\maxi_a$ inside an orbit $\ZZ^n\cdot\maxi_a$ containing a break, we can first translate the whole orbit by $-a$ to the origin. Then rescale in every coordinate direction by $b_i\i$, so that the orbit becomes the standard $\ZZ$-lattice in $\field^n$. In particular, the breaks $g_i,d_i\in (a_i+\ZZ\cdot b_i)$ become points in $\ZZ$ (to be precise, $g_i\up\mapsto \h{g}_i\up= b_i\i(g_i\up-a_i)$, $g_i\low\mapsto \h{g}_i\low= b_i\i(g_i\low-a_i)$). Rescaling and translation are isomorphisms of varieties, so these manipulations are allowed when computing the closure. 
Furthermore, we can compute the closure of $\brac{\maxi}$ over $\QQ$ since $\ol{\brac{\maxi}_\field}=k\otimes_\QQ \ol{\brac{\maxi}_\QQ}$. 
Use the following results from \cite[Section 7.1]{mvdb}:
\begin{prop}\label{prop:closure}
	Consider $\ZZ^n\subset\QQ^n$.
\begin{enumerate}[label=\roman{*}), ref=(\ref{prop:closure}.\roman{*})]
\item\label{lem:techgeo} Given any $\lambda_1,\ldots,\lambda_m \ \in \ (\QQ^n)\ua$,
	there is a unique decomposition of the index set $T=\{1,\ldots,m\}$ into two disjoint parts $I\dot{\cup}J$, such that
	there are $e \in \QQ^n, \ z = (z_1,\ldots,z_m) \in \QQ^m$ with
	\begin{equation*}
	\sum\limits_{i=1}^m z_i \lambda_i \ = \ 0, \;\text{   and   }\;\langle \lambda_i, e\rangle = \lambda_i(e)= \
	\begin{cases} >0, &\text{for } i \in I\\ =0, &\text{for } i \in J \end{cases}\;\text{   and   }\;
	z_i = \ \begin{cases} =0, &\text{for } i \in I\\ >0, &\text{for } i \in J. \end{cases}	
   \end{equation*}
\item\label{prop:geo} Given furthermore $q_1,\ldots,q_m\in\QQ$, define $E=\bigcap_{j\in J} \ker(\lambda_j)$ and
	\begin{eqnarray*}
	C=\{x\in \QQ^n \;|\; \langle \lambda_i,x\rangle=\lambda_i(x) \leq q_i, \; \forall i\in T\},
	&&C\p=\{x\in \QQ^n \; |\; \langle \lambda_j,x\rangle=\lambda_j(x) \leq q_j \; \forall j\in J\},
	\end{eqnarray*}
	then the Zariski closure of $C\cap \ZZ^n$ equals
	$\ol{C\cap \ZZ^n} = C\p \cap (\ZZ^n+E)$
	and	$C\p \cap(\ZZ^n+E)$ is a finite union of translates of $E$.
\item\label{cor:geo}
	For $x\in \ZZ^n$, one has $\ol{(x+C\cap \ZZ^n)\cap (C\cap \ZZ^n)} = \ol{(x+C\cap \ZZ^n)}\cap \ol{(C\cap \ZZ^n)}$.
\end{enumerate}
\end{prop}
This proposition can be applied to the translated, rescaled support of $L(\maxi)$ given by $\ZZ^n\cap C$ with
\begin{align*}
C\ &=\ \{x\in \QQ^n\ |\ \h{g}_i\low+1 \leq x_i \leq \h{g}_i\up\ \text{for all }i\}\\
&=\ \{x\in\QQ^n\ |\ -\varepsilon_i(x) \leq -\h{g}_i\low-1, quad\varepsilon_i(x) \leq \h{g}_i\up,\quad \ 1\leq i\leq n\}\\
&=\ \{x\in\QQ^n\ |\ \lambda_k(x) \leq q_k,\ 1\leq k\leq 2n\}
\end{align*}
where $\varepsilon_i$ denotes the $i$-th coordinate function, $\lambda_k=\varepsilon_k$, $q_k= \h{g}_k\up$ for $1\leq k\leq n$ and $\lambda_k =-\varepsilon_{k-n}$, $q_k=-\h{g}_{k-n}\low-1$ for $n+1\leq k\leq 2n$. Inequalities where $\h{g}_i\up$ or $\h{g}_i\low$ are $\pm\infty$ are dropped.
In our easy situation, we can make the index set $J\subset\{1,\ldots,2n\}$ concrete: 
$$J\ =\ \{i\ |\ \text{ neither }\h{g}_i\up\text{ nor }\h{g}_i\low\ =\ \pm\infty\}$$
(choose eg. $e=(e_k)_k$ with $e_k=e_{k+n}=0$ for those $1\leq k\leq n$ with neither $\h{g}_k\up$ nor $\h{g}_k\low$ are $\pm\infty$, and $e_k=1$ resp. $e_{n+k}=-1$ otherwise. 
Similarly, $z=(z_k)_k$ with $z_k=z_{k+n}=1$ for those $1\leq k\leq n$ with neither $\h{g}_k\up$ nor $\h{g}_k\low$ are $\pm\infty$, and $z_k=0$ otherwise).
We get 
$$\ol{C\cap \ZZ^n}\ =\ \{x\in \QQ^n\ |\ \h{g}_i\low+1 \leq x_i \leq \h{g}_i\up\ \text{for all }i \text{ st. }\h{g}_i\low\text{ and }\h{g}_i\up\neq\pm\infty\} \cap (\ZZ^n+\QQ^{n-J})$$
where we denote $\QQ^{n-J}=\spann{\QQ}\{e_i\ |\ 1\leq i\leq n \text{ and }i\notin J\}$.
Tensor with $\field$ and undo the rescaling and translating, then we get 
\begin{align*}
\ol{\brac{\maxi}}\ &=\ \{x\in \field^n\ |\ g_i\low+b_i \leq x_i \leq g_i\up\ \text{for all }i \in J,\text{ i.e. }g_i\low\text{ and }g_i\up\neq\pm\infty\}\\
& \cap\ (\ZZ^n\cdot b + a + \spann{\field}\{b_i\ |\ 1\leq i\leq n \text{ and }i\notin J\})
\end{align*}
(note here that the inequalities still make sense over an arbitrary field $\field$ because in the $i$-th coordinate for $i\in J$, we work in a lattice).
But because there are only finitely many breaks, there are only finitely many possibilities to choose $g_i\up$ and $g_i\low$ corresponding to a break, as well as for $J\subset \{1,\ldots,n\}$. Therefore there are only finitely many different Zariski closed sets $\ol{\brac{\maxi}}$, so \ref{mvdb5} holds. Finally, \ref{mvdb6} is the consequence of Proposition \ref{cor:geo}.
\begin{bem}
Of course in this easy case the closures can be computed by hands. But this proposition indicates how to deal with (twisted) GWA's where the breaks need no longer be parallel to the coordinate hyperplanes (for twisted GWA's, see \cite{maz-tgwa}).
\end{bem}

\subsection{The proof: The refinement}\label{sec:refine}

Given any primitive ideal $\aideal$, Theorem \ref{thm:mvdb} assigns a simple weight module $L(\maxi)$ such that $\anni{A}(L(\maxi))=\aideal$. Now we show that it is possible to choose $\maxi\p$ to be \emph{highest} weight with $\anni{A}(L(\maxi\p))=\aideal$, under the assumption that none of the $t_i$ is a unit. In that case the tensor factor $A_i$ of $A=A_1\otimes\ldots\otimes A_n$ would be a commutative algebra and not of interest. Once the Theorem gave us $\maxi$, there are two possibilities: 
\begin{itemize}
\item
Either there are breaks $\sigma_i^{\gamma_i\up}(\maxi)$ for $\gamma_i\up>0$ in all coordinate directions $i$. This means that $\sigma^{\gamma\up-\underline{1}}(\maxi)=:\maxi\p$ is a highest weight (where $\underline{1}=(1,\ldots,1)$), and since $\maxi\p$ lies in the support of $L(\maxi)$, we have $L(\maxi)\cong L(\maxi\p)$. Hence $J(\maxi)=J(\maxi\p)$.
\item
Or we have some coordinate $i$ for which $g_i\up=\infty$, so in
\begin{align*}
\ol{\brac{\maxi}}\ &=\ \{x\in \field^n\ |\ g_i\up \geq x_i \geq g_i\low+b_i\ \text{for all }i\in J, \text{ i.e. }g_i\up\text{ and }g_i\low\neq\pm\infty\}\\
& \cap\ (\ZZ^n\cdot b + a + \spann{\field}\{b_i\ |\ 1\leq i\leq n \text{ and }i\notin J\}),
\end{align*}
there is no inequality restricting the coordinate $x_i$ of any element $x\in\ol{\brac{\maxi}}$. In other words, $\ol{\brac{\maxi}}+\field\cdot e_i=\ol{\brac{\maxi}}$. We want to replace $\maxi$ by some other maximal ideal $\maxi\p$ so that their closures are the same, but $L(\maxi\p)$ is a highest weight module. All we need to do is to keep the inequalities and the index set $J$ in the description of $\ol{\brac{\maxi}}$ unchanged.
Replace for this purpose $\maxi=\maxi_a=(T_1-a_1,\ldots,T_n-a_n)$ by any other maximal ideal of the form $(T_1-a_1,\ldots,T_i-z,\ldots,T_n-a_n)$ such that $(T_i-z)$ is a root of $t_i$ (recall that we assumed $t_i\notin\field$). Assume that it is the smallest break in the orbit $\sigma_i^{\ZZ}(T_1-a_1,\ldots,T_i-z,\ldots,T_n-a_n)$. This is possible because $t_i$ has only finitely many roots. Then $\sigma_i(T_1-a_1,\ldots,T_i-z,\ldots,T_n-a_n)=:\maxi\p$ is a highest weight in the $i$-th coordinate direction. Let us check that we preserved the closure $\ol{\brac{\maxi}}=\ol{\brac{\maxi\p}}$: Because we chose the break to be smallest possible, we have $g_i\up=z$ and $g_i\low=-\infty$, and in the computation of the closure the corresponding $i$-th inequality will be dropped. The other coordinate directions are not concerned. Repeating this for all coordinates with $g_i\up=\infty$, we end up with a maximal ideal that is highest weight.
\end{itemize}
Notice that in the last case the two simple modules $L(\maxi)$ and $L(\maxi\p)$ are no longer isomorphic (we even changed the weight lattice), but their annihilators satisfy $J(\maxi) = A\cdot I(\ol{\brac{\maxi}}) = A\cdot I(\ol{\brac{\maxi\p}}) = J(\maxi\p)$, hence the result is the same primitive ideal we started with.

\section{Examples}\label{sec:ex}

In this section, our ground field $\field=\CC$ are the complex numbers.
\subsection{The first Weyl algebra}  
The first Weyl algebra $\CC[x,\partial] = \CC\langle x,\partial\rangle/[\partial,x]=1$ of differential operators on a polynomial ring in one variable can be described as a GWA $A$ of rank one with base ring $R=\CC[T]$, defining element $t=T$ and automorphism $\sigma(T)=T-1$, see \cite[Example 1.2 (1)]{bav}. In particular, since $\sigma$ is a translation with $b:=-1$, it is a GWA of the special form we discuss here.
The defining element $t$ has only one zero, namely $z=0$, hence only one orbit inside $\CC\cong\mspec(\CC[T])$ contains a break, and this is $0+\ZZ\cdot(-1) =\ZZ$.
All modules $M(\maxi_a)$ with $a\notin\ZZ$ are already simple, i.e. $L(\maxi_a)=M(\maxi_a)$, and $\brac{\maxi_a}=a+\ZZ$ is dense in $\CC$, therefore $\anni{A}(L(\maxi_a))= A\cdot I(\ol{\brac{\maxi_a}})\cdot A =(0)$. Instead, concentrate on those $L(\maxi_a)$ with $a\in\ZZ$, eg. $a=2$. The following picture shows the weight lattice of $M(\maxi_2)$:
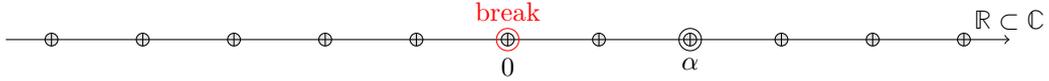
\begin{figure}[H]
		\begin{center}
		\begin{tikzpicture}[scale=1.2]
		
		\coordinate (alph) at (2,0) {};
		
		\draw[->] (-5.5,0) -- (5.5,0) node[above] {$\RR\subset\CC$};
		
		\foreach \x in {-5, ..., 5}
			\draw (\x,-2pt) -- (\x,2pt);
		\draw (0,-3pt) node[anchor=north] {$0$};
		
		\draw[red] (0,0) circle (3.5pt);
		\draw[red] ($(0,0)+(0,+3pt)$) node[anchor=south] {$\text{break}$};
		
		\draw (alph) circle (3.5pt);
		\draw ($(alph)+(0,-3pt)$) node[anchor=north] {$\alpha$};
		
		\foreach \x in {-5, ..., 5}
			\filldraw [fill opacity=0.1] (\x,0) circle (2pt);
					
		\end{tikzpicture}
		\end{center}
\end{figure} 
The action of $X$ and $Y$ on the weight spaces are bijective (gray arrows) except for $M(\maxi_s)_{\maxi_0}$, where the break is: Here $X\cdot M(\maxi_2)_{\maxi_0}=0$.
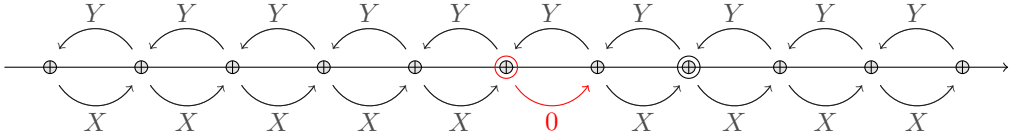
\begin{figure}[H]
		\begin{center}
		\begin{tikzpicture}[scale=1.2]
		
		\coordinate (alph) at (2,0) {};
		
		\draw[->] (-5.5,0) -- (5.5,0);
		
		\foreach \x in {-5, ..., 5}
			\draw (\x,-2pt) -- (\x,2pt);
		\draw[red] (0,0) circle (3.5pt);
		
		\draw (alph) circle (3.5pt);
		
		\foreach \w in {-5, ..., 5}
			\filldraw [fill opacity=0.1] (\w,0) circle (2pt);
		
		\foreach \y in {-5, ...,4}{
		\draw[<-, opacity=0.7] (\y+0.1,0.2) .. controls (\y+0.3,0.5) and (\y+0.7,0.5) .. (\y+0.9,0.2);
		\draw[opacity=0.7] (\y+0.5,0.4) node[anchor=south] {$Y$};}
				
		\foreach \x in {-5, ...,-1}{
		\draw[->, opacity=0.7] (\x+0.1,-0.2) .. controls (\x+0.3,-0.5) and (\x+0.7,-0.5) .. (\x+0.9,-0.2);
		\draw[opacity=0.7] (\x+0.5,-0.4) node[anchor=north] {$X$};}
		
		\draw[->, red] (0.1,-0.2) .. controls (0.3,-0.5) and (0.7,-0.5) .. (0.9,-0.2);
		\draw[red] (0.5,-0.4) node[anchor=north] {$0$};

		\foreach \x in {1, ...,4}{
		\draw[->, opacity=0.7] (\x+0.1,-0.2) .. controls (\x+0.3,-0.5) and (\x+0.7,-0.5) .. (\x+0.9,-0.2);
		\draw[opacity=0.7] (\x+0.5,-0.4) node[anchor=north] {$X$};}
				
		\end{tikzpicture}
		\end{center}
\end{figure}
Thus, $M(\maxi_2)$ has one submodule generated by $M(\maxi_2)_{\maxi_0}$:
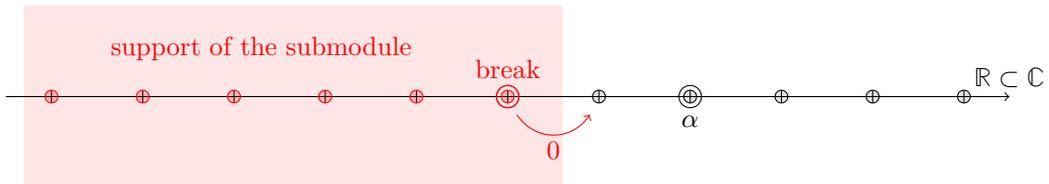
\begin{figure}[H]
		\begin{center}
		\begin{tikzpicture}[scale=1.2]
		
		\coordinate (alph) at (2,0) {};
		
		\fill[fill opacity=0.1, red] (-5.3,-1) -- (0.6,-1) -- (0.6,1) -- (-5.3,1) -- cycle;
		
		\draw[->] (-5.5,0) -- (5.5,0) node[above] {$\RR\subset\CC$};
		
		\foreach \x in {-5, ..., 5}
			\draw (\x,-2pt) -- (\x,2pt);
		
		\draw[red] (0,0) circle (3.5pt);
		\draw[red] ($(-2.7,0)+(0,0.3)$) node[anchor=south] {$\text{support of the submodule}$};
		
		\draw[red] (0,0) circle (3.5pt);
		\draw[red] ($(0,0)+(0,+3pt)$) node[anchor=south] {$\text{break}$};

		\draw[->, red] (0.1,-0.2) .. controls (0.3,-0.5) and (0.7,-0.5) .. (0.9,-0.2);
		\draw[red] (0.5,-0.4) node[anchor=north] {$0$};
		
		\draw (alph) circle (3.5pt);
		\draw ($(alph)+(0,-3pt)$) node[anchor=north] {$\alpha$};
		
		\foreach \x in {1, ..., 5}
			\filldraw [fill opacity=0.1] (\x,0) circle (2pt);
					
		\foreach \x in {-5, ..., 0}
			\filldraw [fill opacity=0.2, red] (\x,0) circle (2pt);
					
		\end{tikzpicture}
		\end{center}

\end{figure}
So we see that also for the two simple weight modules with support $\ZZ_{\leq0}$ resp. $\ZZ_{>0}$, the closure of the support is $\ol{\brac{\maxi_{0}}}=\ol{\brac{\maxi_{1}}}=\CC$ and therefore $\anni{\CC[x,\partial]}(L(\maxi_0))=\anni{\CC[x,\partial]}(L(\maxi_1))=(0)$.
In other words, the only primitive ideal in $\CC[x,\partial]$ is $(0)$, which matches the fact that the Weyl algebra is simple, so there are no nontrivial twosided ideals.
\begin{bem}
Notice that we get a break at $0$, while the computations in \cite[Section 6]{mvdb} correspond to a break at $-1$. This difference can be explained by the choice of $R$. We follow the convention in \cite{bav}, where $R=\field[T]=\field[YX]$, while in \cite{mvdb} $R=\field[T]=\field[XY]$. Since $YX-XY=1$, it follows that 
$$\maxi_0^\text{Bavula} = (YX) = (XY+1) = \maxi_{-1}^\text{MvdB},$$
which explains the `shift by $1$'. The same has to be kept in mind for the $n$-th Weyl algebra.
\end{bem}
\subsection{A rank $1$ example with two breaks}
We stay in the rank $1$ case, we keep the translation $\sigma(T)=T-1$, but we change $t$ to be some other polynomial (these are the `main objects' considered in \cite{bav}). For example, choose $t=(T-3)(T-2)(T+\frac23)(T-(2+i))(T-(4+i))$. Then we have three orbits with breaks: $\ZZ$, $-\frac23+\ZZ$ and $i+\ZZ$.
First we depict how these orbits lie inside the complex plane $\mspec(\CC[T])\cong \CC$ (not to be confused with the following discussion of the rank $2$ case!): 
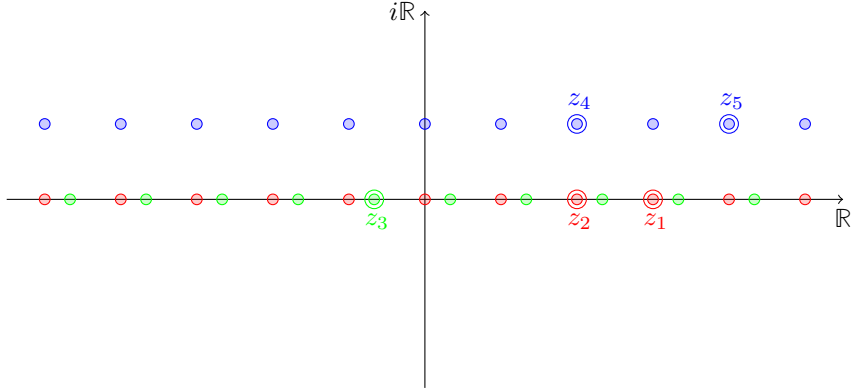
\begin{figure}[H]
		\begin{center}
		\begin{tikzpicture}[scale=1]
		
			\draw[->] (-5.5,0) -- (5.5,0) node[below] {$\RR$};
			\draw[->] (0,-2.5) -- (0,2.5) node[left] {$i\RR$};
			
			\coordinate (z1) at (3,0) {};
			\coordinate (z2) at (2,0) {};
			\coordinate (z3) at (-2/3,0) {};
			\coordinate (z4) at (2,1) {};
			\coordinate (z5) at (4,1) {};
			
			\draw[red] (z1) circle (3.5pt);
			\draw[red] ($(z1)+(1pt,-2pt)$) node[anchor=north] {$z_1$};
			\draw[red] (z2) circle (3.5pt);
			\draw[red] ($(z2)+(1pt,-2pt)$) node[anchor=north] {$z_2$};
			\draw[green] (z3) circle (3.5pt);
			\draw[green] ($(z3)+(1pt,-2pt)$) node[anchor=north] {$z_3$};
			\draw[blue] (z4) circle (3.5pt);
			\draw[blue] ($(z4)+(1pt,2pt)$) node[anchor=south] {$z_4$};
			\draw[blue] (z5) circle (3.5pt);
			\draw[blue] ($(z5)+(1pt,2pt)$) node[anchor=south] {$z_5$};
			
			\foreach \x in {-5, ..., 5}
				\filldraw [fill opacity=0.2, red] (\x,0) circle (2pt);						

			\foreach \x in {-4, ..., 5}
				\filldraw [fill opacity=0.2, green] (\x-2/3,0) circle (2pt);

			\foreach \x in {-5, ..., 5}
				\filldraw [fill opacity=0.2, blue] (\x,1) circle (2pt);						
			\end{tikzpicture}
	\end{center}
	\end{figure} 
Pick the blue orbit, it is the support of eg. $M(\maxi_{0+i})$. Determine its submodules: We have two breaks in the orbit of $0+i$, namely $z_4=2+i$ and $z_5=4+i$. We have observed earlier that for $\alpha>0$,
\begin{align*}
YX^\alpha =0\ &\text{ iff }\ \sigma^{\alpha-1}(t)\in\maxi_{0+i},\\
&\text{ iff }\ (T-(2+i))(T-(4+i))\in\maxi_{\alpha-1+i},\\
&\text{ iff }\ \alpha=3 \text{ or }\alpha=5.
\end{align*} 
$X^3$ and $X^5$ are bijective, while $YX^3$ and $YX^5$ are zero. So there are two submodules, one generated by $X^3$ and the other by $X^5$. This is depicted below, where we shade the support of the two submodules blue.
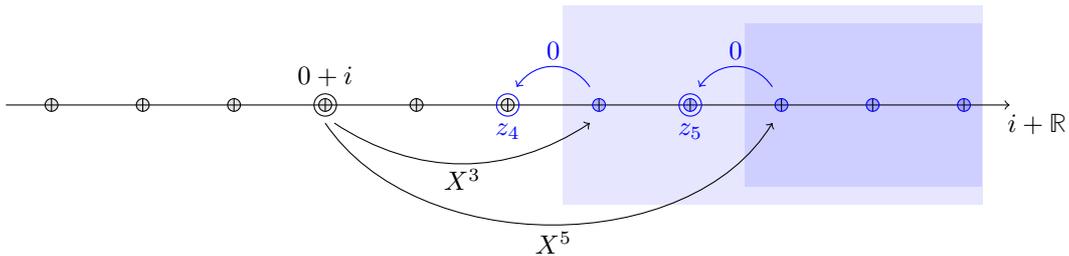
\begin{figure}[H]
		\begin{center}
		\begin{tikzpicture}[scale=1.2]
		
		\coordinate (alph) at (0,0) {};
		\coordinate (b1) at (2,0) {};
		\coordinate (b2) at (4,0) {};
				
		\fill[fill opacity=0.1, blue] (2.6,-1.1) -- (7.2,-1.1) -- (7.2,1.1) -- (2.6,1.1) -- cycle;

		\fill[fill opacity=0.1, blue] (4.6,-0.9) -- (7.2,-0.9) -- (7.2,0.9) -- (4.6,0.9) -- cycle;
				
		\draw[->] (-3.5,0) -- (7.5,0);
		\draw (7.8,0) node[below] {$i+\RR$};
		
		\foreach \x in {-3, ..., 7}
			\draw (\x,-2pt) -- (\x,2pt);

		\draw[blue] (b1) circle (3.5pt);
		\draw[blue] ($(b1)+(0,-3pt)$) node[anchor=north] {$z_4$};

		\draw[blue] (b2) circle (3.5pt);
		\draw[blue] ($(b2)+(0,-3pt)$) node[anchor=north] {$z_5$};

		\draw[<-, blue] (2.1,0.2) .. controls (2.3,0.5) and (2.7,0.5) .. (2.9,0.2);
		\draw[blue] (2.5,0.4) node[anchor=south] {$0$};

		\draw[<-, blue] (4.1,0.2) .. controls (4.3,0.5) and (4.7,0.5) .. (4.9,0.2);
		\draw[blue] (4.5,0.4) node[anchor=south] {$0$};

		\draw[->] (0.1,-0.2) .. controls (1,-0.8) and (2,-0.8) .. (2.9,-0.2);
		\draw (1.5,-0.6) node[anchor=north] {$X^3$};
			
		\draw[->] (0,-0.2) .. controls (1,-1.7) and (4,-1.7) .. (4.9,-0.2);
		\draw (2.5,-1.3) node[anchor=north] {$X^5$};
							
		\draw (alph) circle (3.5pt);
		\draw ($(alph)+(0,3pt)$) node[anchor=south] {$0+i$};
		
		\foreach \x in {-3, ..., 2}
			\filldraw [fill opacity=0.1] (\x,0) circle (2pt);
					
		\foreach \x in {3, ..., 7}
			\filldraw [fill opacity=0.2, blue] (\x,0) circle (2pt);
					
		\end{tikzpicture}
		\end{center}
\end{figure}
Notice that the support of $M(\maxi_{i+3})$ is the same, but the submodule structure is different (still, the subquotients are of course isomorphic):
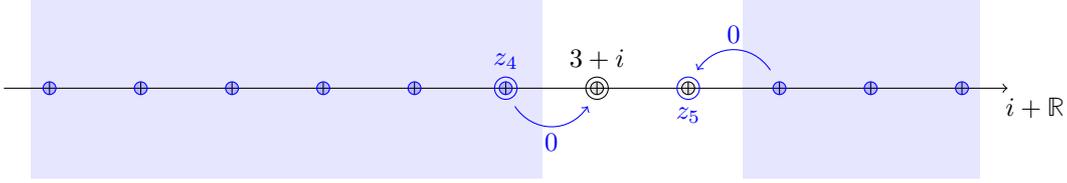
\begin{figure}[H]
		\begin{center}
		\begin{tikzpicture}[scale=1.2]
		
		\coordinate (alph) at (3,0) {};
		\coordinate (b1) at (2,0) {};
		\coordinate (b2) at (4,0) {};
				
		\fill[fill opacity=0.1, blue] (-3.2,-1) -- (2.4,-1) -- (2.4,1) -- (-3.2,1) -- cycle;

		\fill[fill opacity=0.1, blue] (4.6,-1) -- (7.2,-1) -- (7.2,1) -- (4.6,1) -- cycle;
				
		\draw[->] (-3.5,0) -- (7.5,0);
		\draw (7.8,0) node[below] {$i+\RR$};
		
		\foreach \x in {-3, ..., 7}
			\draw (\x,-2pt) -- (\x,2pt);
		
		\draw[blue] (b1) circle (3.5pt);
		\draw[blue] ($(b1)+(0,3pt)$) node[anchor=south] {$z_4$};

		\draw[blue] (b2) circle (3.5pt);
		\draw[blue] ($(b2)+(0,-3pt)$) node[anchor=north] {$z_5$};

		\draw[->, blue] (2.1,-0.2) .. controls (2.3,-0.5) and (2.7,-0.5) .. (2.9,-0.2);
		\draw[blue] (2.5,-0.4) node[anchor=north] {$0$};

		\draw[<-, blue] (4.1,0.2) .. controls (4.3,0.5) and (4.7,0.5) .. (4.9,0.2);
		\draw[blue] (4.5,0.4) node[anchor=south] {$0$};
				
		\draw (alph) circle (3.5pt);
		\draw ($(alph)+(0,3pt)$) node[anchor=south] {$3+i$};
		
		\foreach \x in {3, ..., 4}
			\filldraw [fill opacity=0.1] (\x,0) circle (2pt);
					
		\foreach \x in {-3, ..., 2}
			\filldraw [fill opacity=0.2, blue] (\x,0) circle (2pt);

		\foreach \x in {5, ..., 7}
			\filldraw [fill opacity=0.2, blue] (\x,0) circle (2pt);
					
		\end{tikzpicture}
		\end{center}
\end{figure}
In our notation from above, the lower and upper break for $3+i$ are $g\low=z_4$ and $g\up=z_5$, so the support of $L(\maxi_{3+i})$ is 
$$\brac{\maxi_{3+i}}\ =\ \supp(L(\maxi_{3+i}))\ =\ (i+\ZZ)\ \cap\ \{x\in \CC\ |\ z_5 \geq x > z_4\}\ =\ \{3+i,\ 4+i\}.$$
Since it consists only of two points, it agrees with its closure and hence 
$$\anni{A}(L(\maxi_{3+i}))=A\cdot (\maxi_{3+i}\cap\maxi_{4+i})\cdot A=A\cdot (\maxi_{3+i}\maxi_{4+i})\cdot A.$$
There are up to isomorphism two more simple modules with support in the orbit $i+\ZZ$, namely $L(\maxi_{2+i})$ and $L(\maxi_{5+i})$, both of which have infinite support $i+\ZZ_{\leq 2}$ and $i+\ZZ_{>4}$, resp. The closure of the support is in both cases equal to $\CC$, so the annihilators of both simple modules are $(0)$.
The two other orbits containing breaks can be treated similarly. 
We find only one more nonzero annihilator, namely $\anni{A}(L(\maxi_3))=A\maxi_3 A$, since an orbit needs to contain at least two breaks to allow finite support.

\subsection{A rank $2$ example}
Consider the GWA $A$ with base ring $R=\field[T_1,T_2]$, with automorphisms $\sigma_1(T_1)=T_1-1$, $\sigma_2(T_2)=T_2-\frac{3}{2}$ and with defining elements $t_1=(T_1+2)(T_1-1)$ and $t_2=(T_2+3)(T_2-3)$. Now choose $\maxi=\maxi_{(0,0)}$. The support of $M(\maxi_{(0,0)})$ is given by 
$$\supp(M(\maxi_{(0,0)}))\ =\ (0,0)\ +\ \ZZ\cdot e_1\ +\ \frac{3}{2}\ZZ\cdot e_2,$$
so it contains both breaks $-2$ and $1$ for the first coordinate (corresponding to the maximal ideals $\maxi_{(-2,\alpha_2)}$ and $\maxi_{(1,\alpha_2)}$ for arbitrary $\alpha_2\in \frac{3}{2}\ZZ$) and both breaks $-3$ and $3$ for the second coordinate (corresponding to the maximal ideals $\maxi_{(\alpha_1,-3)}$ and $\maxi_{(\alpha_1,3)}$ for arbitrary $\alpha_1\in \ZZ$). The left picture shows the breaks as (red) hyperplanes in $\field^2$. Since break ideals are those ideals $\maxi$ for which 
$$M_\maxi\xrightarrow{X_i=0}M_{\sigma_i(\maxi)}\quad\text{or}\quad M_\maxi\xleftarrow{Y_i=0}M_{\sigma_i(\maxi)},$$
we furthermore depict $\sigma_i(\text{break in direction }i)$ (light red).
The right picture shows the resulting submodule structure of $M(\maxi_{(0,0)})$:
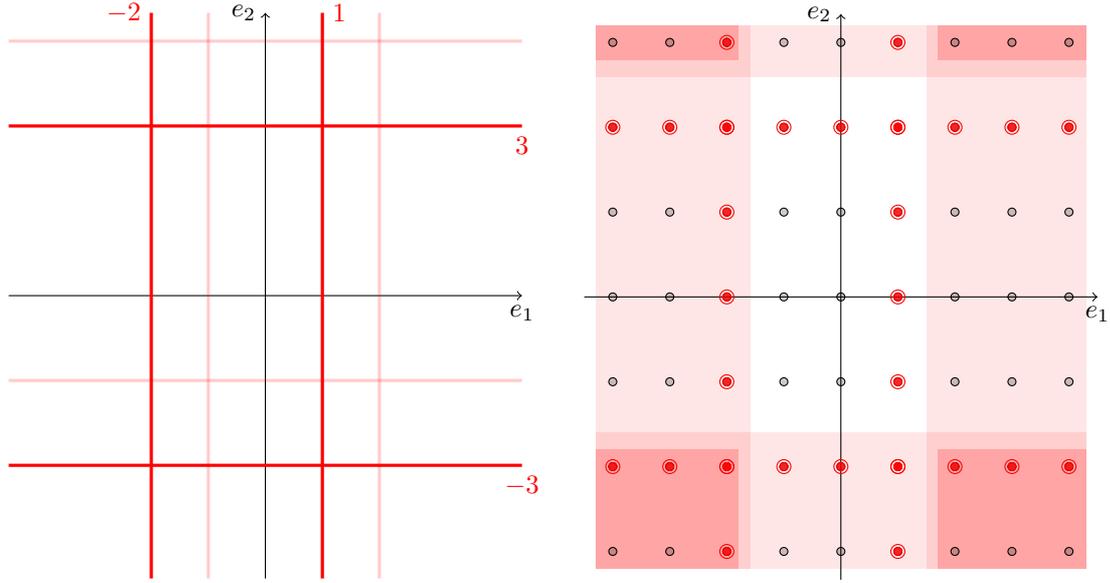
\begin{figure}[H]
		\begin{center}
		\subfigure{
		\begin{tikzpicture}[scale=0.75]
		
			\draw[->] (-4.5,0) -- (4.5,0) node[below] {$e_1$};
			\draw[->] (0,-5) -- (0,5) node[left] {$e_2$};
			
			\draw[red,very thick] (-2,-5) -- (-2,5) node[left] {$-2$};
			\draw[red,very thick,opacity=0.2] (-1,-5) -- (-1,5);

			\draw[red,very thick] (1,-5) -- (1,5) node[right] {$1$};
			\draw[red,very thick,opacity=0.2] (2,-5) -- (2,5);
			
			\draw[red,very thick] (-4.5,3) -- (4.5,3) node[below] {$3$};
			\draw[red,very thick,opacity=0.2] (-4.5,4.5) -- (4.5,4.5);
			
			\draw[red,very thick] (-4.5,-3) -- (4.5,-3)node[below]{$-3$};
			\draw[red,very thick,opacity=0.2] (-4.5,-1.5) -- (4.5,-1.5);
														
			\end{tikzpicture}}\quad
		\subfigure{
		\begin{tikzpicture}[scale=0.75]
			\fill[red, opacity=0.2] (-4.3,-4.8) -- (-1.8,-4.8) -- (-1.8,-2.7) -- (-4.3,-2.7) -- cycle;
			\fill[red, opacity=0.2] (1.7,-4.8) -- (4.3,-4.8) -- (4.3,-2.7) -- (1.7,-2.7) -- cycle;
			\fill[red, opacity=0.1] (-4.3,-4.8) -- (4.3,-4.8) -- (4.3,-2.4) -- (-4.3,-2.4) -- cycle;			
			
			\fill[red, opacity=0.2] (-4.3,4.8) -- (-1.8,4.8) -- (-1.8,4.2) -- (-4.3,4.2) -- cycle;
			\fill[red, opacity=0.2] (1.7,4.8) -- (4.3,4.8) -- (4.3,4.2) -- (1.7,4.2) -- cycle;
			\fill[red, opacity=0.1] (-4.3,4.8) -- (4.3,4.8) -- (4.3,3.9) -- (-4.3,3.9) -- cycle;			

			\fill[red, opacity=0.1] (-4.3,-4.8) -- (-1.6,-4.8) -- (-1.6,4.8) -- (-4.3,4.8) -- cycle;
			\fill[red, opacity=0.1] (1.5,-4.8) -- (4.3,-4.8) -- (4.3,4.8) -- (1.5,4.8) -- cycle;
					
			\draw[->] (-4.5,0) -- (4.5,0) node[below] {$e_1$};
			\draw[->] (0,-5) -- (0,5) node[left] {$e_2$};
			
			\foreach \x in {-4, ..., 4}
			\foreach \y in {-4.5, -3, ..., 4.5}
				\filldraw [fill opacity=0.2, black] (\x,\y) circle (2pt);						

			\foreach \x in {-4, ..., 4}{
				\filldraw [fill opacity=0.8, red] (\x,-3) circle (2pt);
				\draw [red] (\x,-3) circle (3.5pt);
				\filldraw [fill opacity=0.8, red] (\x,3) circle (2pt);
				\draw [red] (\x,3) circle (3.5pt);}			

			\foreach \y in {-4.5, -3, ..., 4.5}{
				\filldraw [fill opacity=0.8, red] (-2,\y) circle (2pt);
				\draw [red] (-2,\y) circle (3.5pt);
				\filldraw [fill opacity=0.8, red] (1,\y) circle (2pt);
				\draw [red] (1,\y) circle (3.5pt);}
											
			\end{tikzpicture}}
	\end{center}
	\end{figure} 
From the break structure, read off the annihilators of the simple modules:
\begin{align*}
\anni{A}(L(\maxi_{(2,\frac92)}))\ &=\ \anni{A}(L(\maxi_{(-2,-3)}))\ =\ \anni{A}(L(\maxi_{(2,-3)}))\ =\ \anni{A}(L(\maxi_{(-2,\frac92)}))\ =\ (0)\\
\anni{A}(L(\maxi_{(0,\frac92)}))\ &=\ \anni{A}(L(\maxi_{(0,-3)}))\ =\ A\cdot ((T_1+1)\cap(T_1)\cap (T_1-1))\cdot A\\
\anni{A}(L(\maxi_{(-2,0)}))\ &=\ \anni{A}(L(\maxi_{(2,0)}))\ =\ A\cdot ((T_2+1)\cap(T_2)\cap(T_2-1)\cap(T_2-2))\cdot A\\
\anni{A}(L(\maxi_{(0,0)}))\ &=\ A\cdot \left(\maxi_{(-1,-\frac{3}{2})}\cap\maxi_{(0,-\frac{3}{2})}\cap \maxi_{(1,-\frac{3}{2})}\right.\\
&\qquad\quad \cap\maxi_{(-1,0)}\cap\maxi_{(0,0)}\cap \maxi_{(1,0)}\\
&\qquad\quad\cap\maxi_{(-1,\frac32)}\cap\maxi_{(0,\frac32)}\cap \maxi_{(1,\frac32)}\\
&\left.\qquad\quad\cap\maxi_{(-1,3)}\cap\maxi_{(0,3)}\cap \maxi_{(1,3)}\right)\cdot A
\end{align*}
There is no further annihilator ideal in $A$ since we considered already all the breaks.
\bibliographystyle{amsalpha}
\addcontentsline{toc}{section}{Literature}
\bibliography{PrimitiveIdealsOfGWA-lit-v1}

\end{document}